\documentclass[a4paper,leqno]{article}
\usepackage{amssymb,amsmath,amsfonts,amsthm,
mathrsfs}
\usepackage{epsfig}

\newtheorem{theorem}{Theorem}[section]

\newtheorem{proposition}[theorem]{Proposition}
\newtheorem{definition}[theorem]{Definition}
\theoremstyle{definition}
\newtheorem{remark}[theorem]{Remark}
\newtheorem{example}[theorem]{Example}

\newcommand{\wt}[1]{\widetilde{#1}}

\newcommand{\Cinf}{\ensuremath{\mathcal{C}^\infty}}
\newcommand{\Cinfc}{\ensuremath{\mathcal{C}^\infty_{\text{c}}}}
\newcommand{\D}{\ensuremath{{\cal D}}}
\renewcommand{\S}{\mathscr{S}}
\newcommand{\E}{\ensuremath{{\cal E}}}

\newcommand{\LL}{\mathcal{L}}

\newcommand{\mb}[1]{\ensuremath{\mathbb{#1}}}
\newcommand{\N}{\mb{N}}

\newcommand{\R}{\mb{R}}
\newcommand{\C}{\mb{C}}

\newcommand{\cl}[1]{\ensuremath{[#1]}}

\newcommand{\G}{\ensuremath{{\cal G}}}

\newcommand{\Gc}{\ensuremath{{\cal G}_\mathrm{c}}}
\newcommand{\Gcinf}{\ensuremath{{\cal G}^\infty_\mathrm{c}}}
\newcommand{\GS}{\G_{{\, }\atop{\hskip-4pt\scriptstyle\S}}\!}
\newcommand{\EM}{\ensuremath{{\cal E}_{M}}}

\newcommand{\Neg}{\mathcal{N}}

\newcommand{\Ginf}{\ensuremath{\G^\infty}}

\newcommand{\GSinf}{\G^\infty_{{\, }\atop{\hskip-3pt\scriptstyle\S}}}

\newcommand{\lara}[1]{\langle #1 \rangle}
\newcommand{\WF}{\mathrm{WF}}
\newcommand{\WFg}{\WF_{\mathrm{g}}}

\newcommand{\supp}{\mathrm{supp}}
\newcommand{\Char}{\ensuremath{\text{Char}}}
\newcommand{\zs}{\setminus 0}
\newcommand{\CO}[1]{\ensuremath{T^*(#1) \zs}}
\newcommand{\ssc}{\mathrm{sc}}

\newcommand{\Ellsc}{\mathrm{Ell}_\ssc}

\newfont{\bigmath}{cmr12 at 13pt}
\newcommand{\PPsi}{\bigmath{\symbol{9}}}
\newcommand{\Oprop}[1]{{\ }_{\mathrm{pr}}^{\hphantom{m}}\text{\PPsi}_{\ssc}^{ #1}}

\newfont{\grecomath}{cmmi12 at 15pt}

\newcommand{\val}{\mathrm{v}} 
\newcommand{\esp}{\mathrm{e}}

\renewcommand{\d}{\ensuremath{\partial}}

\newfont{\bl}{msbm10 scaled \magstep2}

\newcommand{\beq}{\begin{equation}}
\newcommand{\eeq}{\end{equation}}

\newcommand{\isom}{\cong}

\newcommand{\F}{\ensuremath{{\cal F}}}

\newcommand{\notmid}{\mid\kern-0.5em\not\kern0.5em}

\newcommand{\ga}{\gamma}
\newcommand{\Ga}{\Gamma}

\newcommand{\eps}{\varepsilon}

\newcommand{\vphi}{\varphi}

\newcommand{\Om}{\Omega}

\newcommand{\compl}[1]{{#1}^{\mathrm{c}}}

\newcommand{\Sscu}{\underline{\mathcal{S}}_{\,\ssc}}
\newcommand{\Nu}{\underline{\mathcal{N}}}
\newcommand{\Syscu}{{\wt{\underline{\mathcal{S}}}}_{\,\ssc}}

\renewcommand{\val}{\mathrm{v}} 
\renewcommand{\esp}{\mathrm{e}}
\newcommand{\mA}{\mathcal{A}}
\newcommand{\mP}{\mathcal{P}}
\newcommand{\mQ}{\mathcal{Q}}
\newcommand{\mU}{\mathcal{U}}
\newcommand{\mH}{\mathcal{H}}
\newcommand{\M}{\mathcal{M}}
\def\Symp{\widetilde{\mathcal{S}}^m_{\rho,\delta}(\Omega\times\mathbb{R}^p)}

\newcommand{\dslash}{d\hspace{-0.4em}{ }^-\hspace{-0.2em}}

\begin{document}

\title{{\bf On duality theory and pseudodifferential techniques for Colombeau
algebras: generalized delta functionals, kernels and wave front sets}}

\author{Claudia Garetto \\
Dipartimento di Matematica\\ Universit\`a di Torino, Italia\\
\texttt{garettoc@dm.unito.it}\\
\ \\
G\"{u}nther H\"{o}rmann\footnote{Supported by FWF grant P16820-N04}\\
Institut f\"ur Mathematik\\
Universit\"at Wien, Austria\\
\texttt{guenther.hoermann@univie.ac.at}
}
\date{}
\maketitle

\begin{abstract} 
Summarizing basic facts from abstract topological modules over Co\-lom\-be\-au generalized
complex numbers we discuss duality of Colombeau algebras. In particular,
we focus on generalized delta functionals and operator kernels as elements
of dual spaces. A large class of examples is provided by pseudodifferential
operators acting on Colombeau algebras. By a refinement of symbol calculus
we review a new characterization of the wave front set for generalized
functions with applications to microlocal analysis.
\end{abstract}

\section{Topological structures in Colombeau algebras and duality theory}

Apart from some early and inspiring work by Biagioni, Pilipovi\'{c}, Scarpal\'{e}zos 
\cite{Biagioni:90, NPS:98, Scarpalezos:92, Scarpalezos:98, Scarpalezos:00},
topological questions have played a marginal role in the existing Colombeau literature.
However, the recent papers on pseudodifferential operators acting on algebras of generalized
functions \cite{Garetto:04, GGO:03, GH:05} and a preliminary kernel theory introduced in \cite{GGO:03}
motivate a renewed interest in topological issues concerning a class of spaces sufficiently large and general
to model the most common Colombeau algebras. In this section we recall the major points of the theory of 
topological $\wt{\C}$-modules and locally convex topological $\wt{\C}$-modules developed in \cite{Garetto:05a, Garetto:05b}. 
As a topic of particular interest, the foundations of duality theory are provided within this framework, dealing with the $\wt{\C}$-module $\L(\G,\wt{\C})$ of all $\wt{\C}$-linear and continuous functionals on $\G$.
In particular, due to the fact that many algebras of generalized functions can be easily viewed as locally convex topological $\wt{\C}$-modules, all the previous theoretical concepts and results can be applied to the \emph{topological dual of a Colombeau algebra}. We discuss generalized delta functionals as examples of elements of this type of dual space. Among generalized operator kernels those associated with Fourier integral operators have a rich structure, which is of relevance in applications such as microlocal analysis. Section 2 of the paper is devoted to a review of new methods and examples in this direction.

\subsection{Topological and locally convex topological $\wt{\C}$-modules}
We begin by recalling some basic notions. 

The ring $\wt{\C}$ of complex generalized numbers is defined factorizing the algebra $\EM:=\{(u_\eps)_\eps\in\C^{(0,1]}:\, \exists N\in\N\ |u_\eps|=O(\eps^{-N})\, \text{as $\eps\to 0$}\}$ with respect to the ideal $\Neg:=\{(u_\eps)_\eps\in\C^{(0,1]}:\, \forall q\in\N\ |u_\eps|=O(\eps^{q})\, \text{as $\eps\to 0$}\}$ and can be endowed with a structure of a topological ring, making use of the function
\beq
\label{val_M}
\val:\EM\to(-\infty,+\infty]:(u_\eps)_\eps\to\sup\{b\in\R:\ |u_\eps|=O(\eps^b)\, \text{as $\eps\to 0$}\}
\eeq
on $\EM$. The properties of $\val$ enables us to use \eqref{val_M} in defining the \emph{valuation}
\beq
\label{val_C}
\val_{\wt{\C}}(u):=\val((u_\eps)_\eps)
\eeq
of the complex generalized number $u=[(u_\eps)_\eps]$. Let now
\beq
\label{norm_C}
\vert\cdot\vert_\esp := \wt{\C}\to [0,+\infty): u\to \vert u\vert_\esp:=\esp^{-\val_{\wt{\C}}(u)}.
\eeq
The coarsest topology on $\wt{\C}$ such that the map $\vert \cdot\vert_\esp$ is continuous (the ``sharp'' topology on $\wt{\C}$ according to \cite{NPS:98, Scarpalezos:92, Scarpalezos:98, Scarpalezos:00}) is compatible with the ring structure.

A topology $\tau$ on a $\wt{\C}$-module $(\G,+)$ is said to be $\wt{\C}$-linear if the addition $\G\times\G\to\G:(u,v)\to u+v$ and the product $\wt{\C}\times\G\to \G:(\lambda,u)\to\lambda u$ are continuous.
A topological $\wt{\C}$-module $\G$ is a $\wt{\C}$-module with a $\wt{\C}$-linear topology. Clearly $\wt{\C}$ is a topological module on itself.

Our investigation of the topological aspects of a $\wt{\C}$-module is based on the following $\wt{\C}$-linear-adaptation of the notions of absorbent, balanced and convex subsets of a vector space.
\begin{definition}
\label{def_subsets}
A subset $A$ of a $\wt{\C}$-module $\G$ is $\wt{\C}$-absorbent if for all $u\in\G$ there exists $a\in\R$ such that $u\in[(\eps^b)_\eps]A$ for all $b\le a$.\\
$A\subseteq\G$ is $\wt{\C}$-balanced if $\lambda A\subseteq A$ for all $\lambda\in\wt{\C}$ with $\vert\lambda\vert_\esp\le 1$.\\
$A\subseteq\G$ is $\wt{\C}$-convex if $A+A\subseteq A$ and $[(\eps^b)_\eps]A\subseteq A$ for all $b\ge 0$.
\end{definition}
A subset $A$ which is both $\wt{\C}$-balanced and $\wt{\C}$-convex is called absolutely $\wt{\C}$-convex. The $\wt{\C}$-convexity cannot be considered as a generalization of the corresponding concept in vector spaces. In fact the only subset $A$ of $\C$ which is $\wt{\C}$-convex is the trivial set $\{0\}$. In the sequel, we shall simply talk of absorbent, balanced or convex subsets, omitting the prefix $\wt{\C}$, when we deal with $\wt{\C}$-modules. 

As in the classical theory of topological vector spaces every topological $\wt{\C}$-module has a base of absorbent and balanced neighborhoods of the origin. In detail, let $\G$ be a topological $\wt{\C}$-module and $\mathcal{U}$ be a base of neighborhoods of the origin. Then for each $U\in\mathcal{U}$, $(i)$ $U$ is absorbent, $(ii)$ there exists $V\in\mathcal{U}$ with $V+V\subseteq U$, $(iii)$ there exists a balanced neighborhood of the origin $W$ such that $W\subseteq U$ (\cite[Proposition 1.3]{Garetto:05a}). This fact ensures the following characterization: $\G$ is separated if and only if $\cap_{U\in\mathcal{U}}U=\{0\}$. The notion of $\wt{\C}$-convexity introduced above is employed in the definition of \emph{locally convex topological $\wt{\C}$-module}.
\begin{definition}
\label{def_convex_module}
A locally convex topological $\wt{\C}$-module is a topological $\wt{\C}$-module which has a base of $\wt{\C}$-convex neighborhoods of the origin.
\end{definition}
Combining this definition with the properties of bases of neighborhoods of the origin in topological $\wt{\C}$-modules one can easily prove that every locally convex topological $\wt{\C}$-module $\G$ has a base of absolutely convex and absorbent neighborhoods of the origin (\cite[Proposition 1.7]{Garetto:05a}).

We now want to deduce some more information on the topology of $\G$ from the nature of the neighborhoods. Our aim is to find a suitable $\wt{\C}$-version of the Minkowski functional, and therefore a suitable $\wt{\C}$-version of seminorm, which will allow to characterize the topology of a locally convex topological $\wt{\C}$-module $\G$ as the topology determined by the ``$\wt{\C}$''-Minkowski functionals of all absolutely convex and absorbent subsets of $\G$. This is possible by means of the notions of \emph{valuation} and \emph{ultra-pseudo-seminorm}.
\begin{definition}
\label{def_ultra_pseudo}
Let $\G$ be a $\wt{\C}$-module. A \emph{valuation} on $\G$ is a function $\val:\G\to(-\infty,+\infty]$ such that
\begin{trivlist}
\item[(i)] $\val(0)=+\infty$,
\item[(ii)] $\val(\lambda u)\ge \val_{\wt{\C}}(\lambda)+\val(u)$ for all $\lambda\in\wt{\C}$, $u\in\G$,
\item[(ii)'] $\val(\lambda u)= \val_{\wt{\C}}(\lambda)+\val(u)$ for all $\lambda=[(c\eps^a)_\eps]$, $c\in\C$, $a\in\R$, $u\in\G$,
\item[(iii)] $\val(u+v)\ge\min\{\val(u),\val(v)\}$.
\end{trivlist}
An \emph{ultra-pseudo-seminorm} on $\G$ is a function $\mP:\G\to[0,+\infty)$ such that
\begin{trivlist}
\item[(i)] $\mP(0)=0$,
 \item[(ii)] $\mP(\lambda u)\le \vert\lambda\vert_\esp \mP(u)$ for all $\lambda\in\wt{\C}$, $u\in\G$,
\item[(ii)'] $\mP(\lambda u)= \vert\lambda\vert_\esp \mP(u)$ for all $\lambda=[(c\eps^a)_\eps]$, $c\in\C$, $a\in\R$, $u\in\G$,
\item[(iii)] $\mP(u+v)\le\max\{\mP(u),\mP(v)\}$.
\end{trivlist}
\end{definition} 
$\mP(u)=\esp^{-\val(u)}$ is a typical example of an ultra-pseudo-seminorm obtained by means of a valuation on $\G$. An \emph{ultra-pseudo-norm} is an ultra-pseudo-seminorm $\mP$ such that $\mP(u)=0$ implies $u=0$. $\vert\cdot\vert_\esp$ introduced in \eqref{norm_C} is an ultra-pseudo-norm on $\wt{\C}$. In the language of $\wt{\C}$-modules the Minkowski functional of a subset $A$ becomes a ``Minkowski valuation'' $\val_A$.
\begin{proposition}[\cite{Garetto:05a}]
\label{prop_gauge}
Let $A$ be an absolutely convex and absorbent subset of a $\wt{\C}$-module $\G$. Then 
\beq
\label{gauge}
\val_A(u):=\sup\{b\in\R:\, u\in[(\eps^b)_\eps]A\}
\eeq
is a valuation on $\G$. Moreover, for $\mP_A(u):=\esp^{-\val_A(u)}$ and $\eta>0$ the chain of inclusions $
\{u\in\G:\, \mP_A(u)<\eta\}\subseteq[(\eps^{-\log(\eta)})_\eps]A\subseteq\{ u\in\G:\, \mP_A(u)\le\eta\}$ 
holds.
\end{proposition}
We usually call $\mP_A$ the \emph{gauge} of $A$. The properties which characterize an ultra-pseudo-seminorm together with Proposition \ref{prop_gauge} yield the following assertions:
\begin{trivlist}
\item[(i)] the topology induced by the family of ultra-pseudo-seminorms $\{\mP_i\}_{i\in I}$ on $\G$ provides $\G$ with the structure of a locally convex topological $\wt{\C}$-module.
\item[(ii)] In a locally convex topological $\wt{\C}$-module $\G$ the original topology is induced by the family of ultra-pseudo-seminorms $\{\mP_U\}_{U\in\mU}$, where $\mU$ is a base of absolutely convex and absorbent neighborhoods of the origin.
\end{trivlist}
\begin{remark}
An inspection of the neighborhoods of the origin gives some information about ``metrizability'' and ``normability''. More precisely the topology of a separated locally convex topological $\wt{\C}$-module with a countable base of neighborhoods of the origin is induced by a metric invariant under translation and if a separated locally convex topological $\wt{\C}$-module has a bounded neighborhood of the origin then its topology is induced by a ultra-pseudo-norm. By definition a subset $A$ of $\G$ is bounded if and only if it is absorbed by every neighborhood $U$ of the origin, i.e., there exists $a\in\R$ such that $A\subseteq [(\eps^b)_\eps]U$ for all $b\le a$.
\end{remark}
Ultra-pseudo-seminorms are particularly useful in proving the continuity of $\wt{\C}$-linear maps between locally convex topological $\wt{\C}$-modules.
\begin{theorem}[\cite{Garetto:05a}]
Let $(\G,\{\mP_i\}_{i\in I})$ and $(\mH,\{\mQ_j\}_{j\in J})$ be locally convex topological $\wt{\C}$-modules. A $\wt{\C}$-linear map $T:\G\to\mH$ is continuous if and only if it is continuous at the origin if and only if for all $j\in J$ there exists a finite subset $I_0\subseteq I$ and a constant $C>0$ such that for all $u\in\G$
\beq
\label{est_gen_lin}
\mQ_j(Tu) \le C \max_{i\in I_0}\mP_i(u).
\eeq
\end{theorem}
We conclude this subsection focusing our attention on examples of topological $\wt{\C}$-modules from Colombeau's theory of generalized functions, starting with the notion of \emph{strict inductive limit}.
\begin{example}
\label{example_ind}

\bf{Strict inductive limit of locally convex topological $\wt{\C}$-modules}\rm

Consider a $\wt{\C}$-module $\G$ and a sequence $(\G_n)_{n\in\N}$ of submodules of $\G$ such that $\G_n\subseteq\G_{n+1}$ for all $n$ and $\G=\cup_{n\in\N}\G_n$. Assume that $\G_n$ is equipped with a locally convex $\wt{\C}$-linear topology $\tau_n$ such that the topology induced by $\tau_{n+1}$ on $\G_n$ is $\tau_n$. $\G$ endowed with the finest $\wt{\C}$-linear locally convex topology $\tau$ which makes each embedding $\iota_n:\G_n\to\G$ continuous is called \emph{strict inductive limit} of the sequence $(\G_n)_{n\in\N}$ of locally convex topological $\wt{\C}$-modules. This topology is determined by the ultra-pseudo-seminorms $\{\mP_U\}_{U}$ where $U$ is any absolutely convex subset of $\G$ such that $U\cap\G_n$ is a neighborhood of $0$ in $\G_n$ for all $n$. By Proposition 1.21 in \cite{Garetto:05a} we have that the topology $\tau$ on $\G$ induces the original topology $\tau_n$ on each $\G_n$ and that $\G$ is separated if each $\G_n$ is separated. Finally, under the additional assumption of $\G_n$ closed in $\G_{n+1}$ according to $\tau_{n+1}$, one can prove that $\G$ is complete if and only if every $\G_n$ is complete and that $A\subseteq\G$ is bounded if and only if $A$ is contained in some $\G_n$ and bounded there (\cite[Theorems 1.32, 1.26]{Garetto:05a}). We recall that as in the classical theory of inductive limits a $\wt{\C}$-linear map $T:\G\to\wt{\C}$ is continuous if and only if $T_{\vert_{\G_n}}$ is continuous for all $n\in\N$ (\cite[Proposition 1.19]{Garetto:05a}).
\end{example}
\begin{example}

\bf{The $\wt{\C}$-module $\G_E$ of generalized functions based on a locally convex topological vector space $E$}\rm

Let $E$ be a locally convex topological vector space topologized through the family of seminorms $\{p_i\}_{i\in I}$. The elements of  
\beq
\label{defME}
 \M_E := \{(u_\eps)_\eps\in E^{(0,1]}:\, \forall i\in I\,\, \exists N\in\N\quad p_i(u_\eps)=O(\eps^{-N})\, \text{as}\, \eps\to 0\}
\eeq
and
\beq
\label{defNE}
 \Neg_E := \{(u_\eps)_\eps\in E^{(0,1]}:\, \forall i\in I\,\, \forall q\in\N\quad p_i(u_\eps)=O(\eps^{q})\, \text{as}\, \eps\to 0\},
\eeq  
are called $E$-moderate and $E$-negligible, respectively. We define the $\wt{\C}$-modules of \emph{generalized functions based on $E$} as the factor space $\G_E := \M_E / \Neg_E$. Since the growth in $\eps$ of an $E$-moderate net is estimated in terms of any seminorm $p_i$ of $E$, it is natural to introduce the \emph{$p_i$-valuation} of $(u_\eps)_\eps\in\M_E$ as 
\beq
\label{valpi}
\val_{p_i}((u_\eps)_\eps) := \sup\{ b\in\R:\quad p_i(u_\eps)=O(\eps^b)\quad \text{as}\, \eps\to 0\}.
\eeq
It is immediate to see that \eqref{valpi} can be used for defining the \emph{$p_i$-valuation} $\val_{p_i}(u)=\val_{p_i}((u_\eps)_\eps)$ of a generalized function $u=[(u_\eps)_\eps]\in\G_E$. In particular, $\val_{p_i}$ is a valuation in the sense of Definition \ref{def_ultra_pseudo} and thus $\mP_i(u):=\esp^{-\val_{p_i}(u)}$ is an ultra-pseudo-seminorm on the $\wt{\C}$-module $\G_E$. $\G_E$ endowed with the topology of the ultra-pseudo-seminorms $\{\mP_i\}_{i\in I}$ is a locally convex topological $\wt{\C}$-module. Following \cite{NPS:98, Scarpalezos:92, Scarpalezos:98, Scarpalezos:00} we use the adjective ``sharp'' for the topology induced by the ultra-pseudo-seminorms $\{\mP_i\}_{i\in I}$. The sharp topology on $\G_E$, here denoted by $\tau_\sharp$ is independent of the choice of the family of seminorms which determines the original locally convex topology on $E$. $(\G_E,\tau_\sharp)$ is a separated locally convex topological $\wt{\C}$-module and it is complete when $E$ has a countable base of neighborhoods of the origin (\cite[Propositions 3.2, 3.4]{Garetto:05a}).
\end{example}
\begin{example}

\bf{Colombeau algebras obtained as $\wt{\C}$-modules $\G_E$}\rm

Particular choices of $E$ give several known Colombeau algebras of generalized functions as $\wt{\C}$-modules $\G_E$ and the corresponding sharp topologies. This is of course the case for $E=\C$ and $\G_E=\wt{\C}$ which is an ultra-pseudo-normed $\wt{\C}$-module.\\ 
Consider now an open subset $\Om$ of $\R^n$. $E=\E(\Om)$, i.e. the space $\Cinf(\Om)$ topologized through the family of seminorms $p_{K_i,j}(f)=\sup_{x\in K_i, |\alpha|\le j}|\partial^\alpha f(x)|$, where $K_0\subset K_1\subset.... K_i\subset...$ is a countable and exhausting sequence of compact subsets of $\Om$, provides $\G_E=\G(\Om)$. $\G(\Om)$ endowed with the sharp topology determined by $\{\mP_{K_i,j}\}_{i\in\N, j\in\N}$ is a Fr\'echet $\wt{\C}$-module (i.e., metrizable and complete). Another example of a Fr\'echet $\wt{\C}$-module is given by $\G_E$ when $E$ is $\S(\R^n)$. In this way we construct the algebra $\GS(\R^n):=\G_{\S(\R^n)}$ (see \cite{Garetto:04, GGO:03}).
\end{example}
\begin{example}
\bf{The Colombeau algebra $\Gc(\Om)$ of compactly supported generalized functions}\rm

For $K\Subset\Omega$ we denote by $\G_K(\Om)$ the space of all generalized functions in $\G(\Om)$ with support contained in $K$. Note that $\G_K(\Om)$ is contained in $\G_{\D_{K'}(\Om)}$ for all compact subsets $K'$ of $\Om$ such that $K\subseteq {\rm{int}}(K')$, where $\D_{K'}(\Om)$ is the space of all smooth functions $f$ with $\supp\, f\subseteq K'$. We recall that in $\G(\Om)$ the $p_{K,n}$-valuation where $p_{K,n}(f)=\sup_{x\in K, |\alpha|\le n}|\partial^\alpha f(x)|$ is obtained as the valuation of the complex generalized number $\sup_{x\in K, |\alpha|\le n}|\partial^\alpha u(x)|:=(\sup_{x\in K, |\alpha|\le n}|\partial^\alpha u_\eps(x)|)_\eps +\Neg$. Hence for $K,K'\Subset\Om$, $K\subseteq\text{int}(K')$,
\beq
\label{valGK}
\val_{K,n}(u)=\val_{p_{K',n}}(u)
\eeq
is a valuation on $\G_K(\Om)$. More precisely \eqref{valGK} does not depend on $K'$. One can easily check that 
$\G_K(\Om)$, with the topology induced by the ultra-pseudo-seminorms $\{\mP_{\G_K(\Om),n}(u):=\esp^{-\val_{K,n}(u)}\}_{n\in\N}$, is a locally convex topological $\wt{\C}$-module and by construction its topology coincides with the topology induced by any $\G_{\D_{K'}(\Om)}$ with $K\subseteq \text{int}(K')$. In particular, $(\G_{K}(\Om),\{\mP_{\G_K(\Om),n}\}_{n\in\N})$ is a Fr\'echet $\wt{\C}$-module. Let now  $(K_n)_{n\in\N}$ be an exhausting sequence of compact subsets of $\Omega$ such that $K_n\subseteq K_{n+1}$. Clearly $\G_c(\Om)=\cup_{n\in\N}\G_{K_n}(\Om)$ and the assumptions of Example \ref{example_ind} are satisfied by $\G_n=\G_{K_n}(\Om)$. Therefore $\G_c(\Om)$ endowed with the strict inductive limit topology of the sequence $(\G_{K_n}(\Om))_{n}$ is a separated and complete locally convex topological $\wt{\C}$-module. 
\end{example}
\begin{example}

\bf{Algebras of regular generalized functions}\rm

For any locally convex topological vector space $(E,\{p_i\}_{i\in I})$ the set
\beq
\label{MinfE}
\M^\infty_E:=\{(u_\eps)_\eps\in E^{(0,1]}:\ \exists N\in\N\, \forall i\in I\quad p_i(u_\eps)=O(\eps^{-N})\ \text{as}\ \eps\to 0\}
\eeq
is a subspace of the set $\M_E$ of $E$-moderate nets. Therefore the corresponding factor space $\G^\infty_E:=\M^\infty_E /\Neg_E$ is a subspace of $\G_E$ whose elements are called \emph{regular generalized functions based on $E$}. The moderateness properties of $\M_E^\infty$ allow us to define the valuation $\val^\infty_E:\M_E^\infty\to (-\infty,+\infty]$ as 
\[
\val^\infty_E((u_\eps)_\eps)= \sup\{b\in\R:\ \forall i\in I\ p_i(u_\eps)=O(\eps^b)\ \text{as}\ \eps\to 0\}
\]
which can be obviously extended to $\G^\infty_E$. This yields the existence of the ultra-pseudo-norm $\mP_E^\infty(u):=\esp^{-\val^\infty_E(u)}$ on $\G_E^\infty$ which determines a topology finer than the one induced by $\G_E$. One can prove that when $E$ has a countable base of neighborhoods of the origin then the associated space $\G^\infty_E$ of regular generalized functions is a complete and ultra-pseudo-normed $\wt{\C}$-module. A concrete example of $\G^\infty_E$ is given by the Colombeau algebra of $\S$-regular generalized functions $\GSinf(\R^n)$ (see \cite{Garetto:04, GGO:03}), whose definition is precisely $\G^\infty_E$ with $E=\S(\R^n)$.

The Colombeau algebra $\Ginf(\Om)$ of regular generalized functions can be seen as the intersection $\displaystyle\cap_{K\Subset\Omega}\Ginf(K)$, where $\Ginf(K)$ is the space of all $u\in\G(\Om)$ such that there exists a representative $(u_\eps)_\eps$ satisfying the condition: 
\[
\exists N\in\N\ \forall\alpha\in\N^n \sup_{x\in K}|\partial^\alpha u_\eps(x)|=O(\eps^{-N})\quad \text{as}\ \eps\to 0.
\]
On $\Ginf(K)$ we define the ultra-pseudo-seminorm $\mP_{\Ginf(K)}$ via the valuation $
\val_{\Ginf(K)}(u)=\sup\{b\in\R:\, \forall\alpha\in\N^n\ \sup_{x\in K}|\partial^\alpha u_\eps(x)|=O(\eps^b)\}$.
It follows that $\Ginf(\Om)$ is a Frech\'et $\wt{\C}$-module if topologized trough the family of ultra-pseudo-seminorms $\mP_{\Ginf(K_n)}$, where $\{K_n\}_{n\in\N}$ is an arbitrary exhausting sequence of compact subsets of $\Om$.

Finally, a strict inductive limit procedure equips $\Gcinf(\Om)=\cup_{n\in\N}\Ginf_{K_n}(\Om)$ with a complete and separated locally convex topology. Each $\Ginf_{K_n}(\Om):=\G_K(\Om)\cap\Ginf(\Om)$ is a complete ultra-pseudo-normed $\wt{\C}$-module with respect to $\mP_{\Ginf_{K_n}(\Om)}(u)$\\
$=\esp^{-\val^\infty_{K_n}(u)}$, where $\val^\infty_{K_n}(u)=\sup\{b\in\R:\, \forall\alpha\in\N^n\, \sup_{x\in K'_n}|\partial^\alpha u_\eps(x)|=O(\eps^b)\}$ and $K_n\subseteq\text{int}(K'_n)\subseteq K'_n\Subset\Om$.
\end{example}

\subsection{Duality theory for topological $\wt{\C}$-modules}
This subsection is devoted to the dual of a topological $\wt{\C}$-module $\G$ i.e. the $\wt{\C}$-module $\LL(\G,\wt{\C})$ of all $\wt{\C}$-linear and continuous maps on $\G$ with values in $\wt{\C}$. Since $\G$ with $\LL(\G,\wt{\C})$ forms a pairing, the topological dual of $\G$ can be endowed with the weak topology $\sigma(\LL(\G,\wt{\C}),\G)$ that is the coarsest topology such that each map $\lara{u,\cdot}:\LL(\G,\wt{\C})\to\wt{\C}:T\to T(u)$ is continuous for $u$ varying in $\G$. This is a separated locally convex $\wt{\C}$-linear topology determined by the family of ultra-pseudo-seminorms $\mP_u(T)=|T(u)|_\esp$. 

Before introducing other topologies on the dual $\LL(\G,\wt{\C})$ we recall that the polar of a subset $A$ of $\G$ is the set $A^\circ$ of all $T\in\LL(\G,\wt{\C})$ such that $|T(u)|_\eps\le 1$ for all $u\in A$. As proved in \cite[Proposition 2.4]{Garetto:05a} $A^\circ$ is a balanced and convex subset of $\LL(\G,\wt{\C})$ closed with respect to $\sigma(\LL(\G,\wt{\C}),\G)$. Moreover $A^\circ$ is absorbent if and only if $A$ is bounded in $(\G,\sigma(\G,\LL(\G,\wt{\C})))$. It follows that for any $\sigma(\G,\LL(\G,\wt{\C}))$-bounded subset of $\G$ we can define the gauge $\mP_{A^\circ}$ according to Proposition \ref{prop_gauge}.
\begin{definition}
\label{def_strong_top}
Let $\G$ be a topological $\wt{\C}$-module. We call strong topology the topology $\beta(\LL(\G,\wt{\C}),\G)$ determined by the ultra-pseudo-seminorms $\{\mP_{A^\circ}\}_{A\in\mA}$ where $\mA$ is the collection of all $\sigma(\G,\LL(\G,\wt{\C}))$-bounded subsets of $\G$.
Restricting $\mA$ to the family of all bounded subset of $\G$ we define the topology $\beta_b(\LL(\G,\wt{\C}),\G)$ of uniform convergence on bounded subsets of $\G$.
\end{definition}
It is clear that $\sigma(\LL(\G,\wt{\C}),\G)\preceq\beta_b(\LL(\G,\wt{\C}),\G)\preceq\beta(\LL(\G,\wt{\C}),\G)$,
where $\preceq$ stands for ``is coarser than''. Proposition 3.26 in \cite{Garetto:05a} proves that $\beta_b(\LL(\G,\wt{\C}),\G)$ and $\beta(\LL(\G,\wt{\C}),\G)$ coincide when $\G$ is a $\wt{\C}$-module of generalized functions based on a normed space $E$. We now focus on some special families of topological $\wt{\C}$-modules. In order to state the following definition we recall that a subset $S$ of $\G$  is said to be bornivourus if it absorbs every bounded subset of $\G$. Finally $S$ is a barrel of $\G$ if it is absorbent, balanced, convex and closed. 
\begin{definition}
A locally convex topological $\wt{\C}$-module $\G$ is \emph{bornological} if every balanced, convex and bornivorous subset of $\G$ is a neighborhood of the origin. A locally convex topological $\wt{\C}$-module is \emph{barrelled} if every barrel is a neighborhood of the origin.
\end{definition} 
The Frech\'et $\wt{\C}$-modules are examples of bornological and barrelled $\wt{\C}$-modules as well as the strict inductive limit of Frech\'et $\wt{\C}$-modules (see \cite[Propositions 2.9, 2.14, 2.15]{Garetto:05a}). Moreover by Proposition 2.10 in \cite{Garetto:05a} we have that the dual of a bornological $\wt{\C}$-module is complete with respect to both the topologies $\beta(\LL(\G,\wt{\C}),\G)$ and $\beta_b(\LL(\G,\wt{\C},\G)$.

\subsection{Topological dual of a Colombeau algebra: generalized delta functionals and operator kernels} 

The theoretical background provided by the previous subsection can be applied to the topological duals of the Colombeau algebras $\Gc(\Om)$, $\G(\Om)$ and $\GS(\R^n)$. Since $\G(\Om)$ and $\GS(\R^n)$ are Frech\'et $\wt{\C}$-modules and $\Gc(\Om)$ is a strict inductive limit of Frech\'et $\wt{\C}$-modules we know that the duals $\LL(\Gc(\Om),\wt{\C})$, $\LL(\G(\Om),\wt{\C})$ and $\LL(\GS(\R^n),\wt{\C})$ are complete locally convex topological $\wt{\C}$-modules when equipped with the strong topology or with the topology of uniform convergence on bo\-un\-ded subsets. 
As an analogy with distribution theory we recall that $\LL(\Gc(\Om),\wt{\C})$ is a sheaf on $\Om$ and that $\LL(\G(\Om),\wt{\C})$ can be identified with the set of all functionals in $\LL(\Gc(\Om),\wt{\C})$ having compact support (\cite[Theorem 1.2]{Garetto:05b}).

Point value theory for generalized functions \cite{GKOS:01, OK:99} and a useful characterization of the ideals of negligible nets which occur in the definition of some Colombeau algebras as $\G(\Om)$, $\GS(\R^n)$, $\G_{p,p}(\R^n):=\G_{W^{\infty,p}(\R^n)}$ (see \cite{Garetto:05b, GKOS:01}) are the tools employed in proving the following embedding theorem. The Colombeau algebras considered in the sequel are endowed with the topologies discussed in Subsection 1.1 and the corresponding duals are equipped with the topology $\beta_b$ of uniform convergence on bounded subsets. It is clear that the identity between representatives defines the $\wt{\C}$-linear continuous embeddings $\Gc(\Om)\subseteq\G(\Om)$, $\Ginf(\Om)\subseteq\G(\Om)$, $\Gcinf(\Om)\subseteq\Gc(\Om)$ and $\GSinf(\R^n)\subseteq\GS(\R^n)$.
\begin{theorem}[\cite{Garetto:05b}]
\label{th_embedding}
By integration we obtain $\wt{\C}$-linear continuous embeddings, $u\to\big(v\to\int_\Om u(y)v(y)\, dy\big)$, of
\begin{itemize}
\item[(i)] $\G(\Om)$ into $\LL(\Gc(\Om),\wt{\C})$,
\item[(ii)] $\Gc(\Om)$ into $\LL(\G(\Om),\wt{\C})$,
\item[(iii)] $\GS(\R^n)$ into $\LL(\GS(\R^n),\wt{\C})$.
\end{itemize}
\end{theorem}
Theorem \ref{th_embedding} yields the following chains of continuous embeddings $\Ginf(\Om)\subseteq\G(\Om)\subseteq\LL(\Gc(\Om),\wt{\C})$, $\Gc^\infty(\Om)\subseteq\Gc(\Om)\subseteq\LL(\G(\Om),\wt{\C})$, $\GSinf(\R^n)\subseteq\GS(\R^n)\subseteq\LL(\GS(\R^n),\wt{\C})$, which play an important role in regularity theory of differential and pseu\-do\-dif\-fe\-ren\-tial operators with generalized symbols (c.f. \cite{Garetto:04th, GGO:03}). Concluding we discuss some interesting examples of functionals in $\LL(\Gc(\Om),\wt{\C})$, $\LL(\G(\Om),\wt{\C})$ and $\LL(\GS(\R^n),\wt{\C})$. They are mainly provided by point value theory and kernel theory for Fourier integral operators with generalized symbols.
 

\begin{example}

\bf{The generalized delta functional $\delta_{\wt{x}}$}\rm

Let $\wt{x}$ be a generalized point in $\wt{\Om}_{\rm{c}}$. We can define a $\wt{\C}$-linear map $\delta_{\wt{x}}:\G(\Om)\to\wt{\C}$ associating with each $u\in\G(\Om)$ its point value $u(\wt{x})$ at $\wt{x}$. It is clear that $\delta_{\wt{x}}$ belongs to $\LL(\G(\Om),\wt{\C})$ since  $\vert\delta_{\wt{x}}(u)\vert_\esp\le \esp^{-\val(\sup_{x\in K}|u_\eps(x)|)}$, where $(x_\eps)_\eps$ is a representative of $\wt{x}$ contained in a compact set $K$ for $\eps$ small enough.

Since every $u\in\Gc(\Om)$ is a generalized function in some $\G_{\D_{K'}(\Om)}$ it is meaningful to define the point value $u(\wt{x})$ even when $\wt{x}$ is in $\wt{\Om}\setminus\wt{\Om}_{\rm{c}}$. In this case
$\delta_{\wt{x}}$ is a $\wt{\C}$-linear map from $\Gc(\Om)$ into $\wt{\C}$ and by $\vert\delta_{\wt{x}}(u)\vert_\esp\le\mP_{K,0}(u)$ on $\G_K(\Om)$ it is continuous. Finally, as proved in \cite{Garetto:05b}, a point value theory with generalized points in $\wt{\R^n}$ can be formulated for generalized functions in $\GS(\R^n)$. One easily sees that for $\wt{x}\in\wt{\R^n}$ the functional $\delta_{\wt{x}}:\GS(\R^n)\to\wt{\C}$ belongs to $\LL(\GS(\R^n),\wt{\C})$.
\end{example}
\begin{remark}

\bf{Properties of $\delta_{\wt{x}}$}\rm

Let $\wt{x}=[(x_\eps)_\eps]\in\wt{\Om}$. Defining $\supp\, \wt{x}$, the support $\supp\, \wt{x}$ of $\wt{x}$, to be the complement of the set
\beq
\label{supp_wtx}
\left\{ x_0\in\Om:\ \begin{array}{cc}\exists V(x_0)\ \text{neighborhood of $x_0 $}\ \ \exists \eta\in(0,1]\ \ \forall\eps\in(0,\eta]\\[0.2cm]
x_\eps\not\in V(x_0)
\end{array}
\right\},  
\eeq
in $\Om$, one can prove that $\supp\, \delta_{\wt{x}}=\supp\, \wt{x}$. Note that by \eqref{supp_wtx} the support of $\wt{x}\in\wt{\Om}$ is the set of all accumulation points of a representing net $(x_\eps)_\eps$. As a consequence and differently from the distributional case, $\delta_{\wt{x}}$ may have support empty or unbounded. Another difference to distribution theory consists in the existence of elements of $\LL(\G(\Om),\wt{\C})$ having support $\{0\}$ which are not $\wt{\C}$-linear combinations of $\delta_0$ and its derivatives. As an example take $[(\eps)_\eps]\in\wt{\R}$ and $\delta_{[(\eps)_\eps]}\in\LL(\G(\R),\wt{\C})$. By \eqref{supp_wtx} $\supp\, \delta_{[(\eps)_\eps]}=\{0\}$ but for any $m\in\N$ and any choice of $c_i\in\wt{\C}$ the equality $\delta_{[(\eps)_\eps]}=\sum_{i=0}^m c_i\delta_{0}^{(i)}$ does not hold.
\end{remark}
The following theorem proves that the restriction of the generalized delta functional $\delta_{\wt{x}}$ to the set of regular generalized functions is an integral operator. This implies an integral representation of regular generalized functions as in \cite[Theorems 2.3, Propositions 2.6]{Garetto:05b}. 
\begin{theorem}
\label{delta_theorem}
\begin{itemize}
\item[{\ }]
\item[(i)] For all $\wt{x}\in\wt{\Om}_{\rm{c}}$ there exists $v\in\Gc(\Om)$ such that for all $u\in\Ginf(\Om)$
\beq
\label{int_representation_formula}
u(\wt{x})=\int_\Om v(y)u(y)\, dy.
\eeq
\item[(ii)] $(i)$ holds with $\wt{\Om}$, $\G(\Om)$ and $\Gcinf(\Om)$ in place of $\wt{\Om}_{\rm{c}}$, $\Gc(\Om)$ and $\Ginf(\Om)$ respectively.
\item[(iii)] For all $\wt{x}\in\wt{\R^n}$ there exists $v\in\GS(\R^n)$ such that for all $u\in\GSinf(\R^n)$
\beq
\label{int_representation_formula_2}
u(\wt{x})=\int_{\R^n} v(y)u(y)\, dy.
\eeq
\end{itemize}
\end{theorem}
The convolution with a mollifier in $\S(\R^n)$ is the basic idea in proving Theorem \ref{delta_theorem}. Indeed, in $(iii)$ we have that $v$ is given by $[(\varphi_\eps(x_\eps-\cdot))_\eps]\in\GS(\R^n)$ where $(x_\eps)_\eps$ is any representative of $\wt{x}$. In general \eqref{int_representation_formula} and \eqref{int_representation_formula_2} do not hold when $u$ is not regular. See \cite[Remarks 2.2, 2.5]{Garetto:05b} for details. 
\begin{example}

\bf{Generalized kernels}\rm

The $\wt{\C}$-linear functional on $\Gc(\Om)$ defined by the (generalized) integral $$\int_\Om k(y)u(y) dy$$ is called integral operator with kernel $k\in\G(\Om)$. It is an element of $\LL(\Gc(\Om),\wt{\C})$, since with $u\in\G_K(\Om)$, $K\subseteq{\rm{int}}(K')\subseteq K'\Subset\Om$, the estimate
\[
\biggl|\int_\Om k_\eps(y)u_\eps(y)\, dy\biggr|\le c\sup_{y\in K'}|k_\eps(y)| \sup_{y\in K'}|u_\eps(y)|
\]
yields
\[
\biggl|\int_\Om k(y)u(y)\, dy\biggr|_\esp\le C\,\mP_{\G_K(\Om),0}(u),
\]
where $C=\esp^{-\val(\sup_{y\in K'}|k_\eps(y)|)}$.\\
Analogously we prove that $u\to\int_\Om k(y)u(y)\, dy$ belongs to $\LL(\G(\Om),\wt{\C})$ when $k\in\Gc(\Om)$ and that the functional $u\to\int_{\R^n}k(y)u(y)dy$ is an element of $\LL(\GS(\R^n),\wt{\C})$ for $k\in\GS(\R^n)$.

Let us now consider a phase function $\phi$ on $\Om\times\R^p$ and a generalized symbol $a=[(a_\eps)_\eps]\in\Symp$ as defined in \cite[Definition 3.3]{GGO:03} with $\rho>0$ and $\delta<1$. The generalized oscillatory integral 
\[
I_{\phi}(au):=\int_{\Omega\times\mathbb{R}^p}\hskip-10pt e^{i\phi(y,\xi)}a(y,\xi)u(y)\hskip2pt dy\,d\xi:=\biggl[\biggl(\int_{\Om\times\R^p}\hskip-10pt e^{i\phi(y,\xi)}a_\eps(y,\xi)u_\eps(y)\hskip2pt dy\,d\xi\biggr)_\eps\biggr]
\]
is an example of a functional in $\LL(\Gc(\Om),\wt{\C})$. In fact, denoting the symbol seminorm $\sup_{y\in K', \xi\in\R^p}\sup_{|\alpha+\beta|\le j}\lara{\xi}^{-m+\rho|\alpha|-\delta|\beta|}|\partial^\alpha_\xi\partial^\beta_y a_\eps(y,\xi)|$ by $|a_\eps|^{(m)}_{K',j,\rho,\delta}$ and working on $I_{\phi}(au)$ at the level of representatives we have that  
\[
\biggl|\int_{\Om\times\R^p}\hskip-10pt e^{i\phi(y,\xi)}a_\eps(y,\xi)u_\eps(y)\, dy\, d\xi\biggr|\le c_{j,K} |a_\eps|^{(m)}_{K',j,\rho,\delta}\sup_{y\in K',|\gamma|\le j}|\partial^\gamma u_\eps(y)| 
\]
is va\-lid for $u\in\G_K(\Om)$, $m-js<-p$ and for some constant $c_{j,K}$ depending on $j$ and $K$. 

When $\phi(x,y,\xi)$ is a phase function on $\Om'\times\Om\times\R^p$ and, in addition, a phase function in $(y,\xi)$ for all $x$, the map $A:\Gc(\Om)\to\G(\Om')$ given by
\begin{multline*}
Au(x)=\int_{\Om\times\R^p}e^{i\phi(x,y,\xi)}a(x,y,\xi)u(y)\, dy\, \dslash\xi\\
:=\biggl[\biggl(\int_{\Om\times\R^p}e^{i\phi(x,y,\xi)}a_\eps(x,y,\xi)u_\eps(y)\, dy\, \dslash\xi\biggr)_\eps\biggr]
\end{multline*}
is called Fourier integral operator with generalized amplitude $a\in\widetilde{S}^m_{\rho,\delta}(\Om'\times\Om\times\R^p)$ 
(\cite[Proposition 3.10]{GGO:03}). The kernel of $A$ is the $\wt{\C}$-linear functional on $\Gc(\Om'\times\Om)$ determined by 
\[
k_A(u):=\int_{\Om'} A(u(x,\cdot))\, dx.
\]
Since $k_A$ can be expressed by the oscillatory integral
\[
\int_{\Om'\times\Om\times\R^p}e^{i\phi(x,y,\xi)}a(x,y,\xi)u(x,y)\, dy\, dx\, \dslash\xi
\]
in $dx$, $dy$, $d\xi$, from the previous reasoning it follows that $k_A$ belongs to the dual $\LL(\Gc(\Om'\times\Om),\wt{\C})$. We recall that when $\Om'=\Om$ and $\phi(x,y,\xi)=(x-y)\xi$ we obtain the kernel of a generalized pseudodifferential operator.

\end{example}

\section{Pseudodifferential techniques and microlocal analysis}

\subsection{The wave front set}

Let $u$ be a distribution or generalized function on $\Om$. In both
cases, the basic idea of the \emph{wave front set} of $u$, $\WF(u)$, 
can be sketched as follows: it consists of pairs 
$(x;\xi)$ in $\Omega
\times \R^d \isom T^*\Om$, where $x_0$ marks the location of
singularities in $\Om$, 
whereas $\xi_0 \not= 0$ in $\R^d$ represents the directions of
high-frequency contributions in their Fourier spectrum. The
focusing procedure for this ``localization
with attached directional spectral analysis'' can be based on spatial
cut-off functions $\vphi\in\D(\Om)$ and cones
$\Ga\subseteq \R^d \zs$, where rapid decrease of the corresponding
Fourier transform $\F(\vphi u)$ is to be tested. To 
be more precise, the \emph{microlocal regularity test} for $u =
\cl{(u_\eps)_\eps}$ 
at $(x_0;\xi_0) \in \CO{\Om}$ (cotangent space deprived of its zero
section) at scale $0 < \ga_\eps \nearrow \infty$ ($\eps \to 0$)
requires the following: find an open neighborhood $U$ of 
$x_0$ and an open cone $\xi_0 \in \Ga \subseteq \R^d \zs$ such that for all
$\vphi\in\D(U)$ the generalized Fourier transform $\F (\vphi u)$ is
$\G$-rapidly decreasing with scale $\ga$ in $\Ga$, 
i.e., $\exists N$ $\forall l$ $\exists \eps_0$ such that
$$
  |\F(\vphi u_\eps)(\xi)| \leq \ga_\eps^N (1 + |\xi|)^{-l} \qquad
  \forall \xi\in\Ga, \forall \eps \in (0,\eps_0).
$$
The complement of the set of regular points, denoted by
$\WF_g^\ga(u)$, is the \emph{generalized wave front set} of
$u$ (at scale $\ga$). If $v\in\D'(\Om)$ then we have full consistency
with the distribution theoretic notion in the sense that $\WF_g^\ga(\iota(v)) =
\WF(v)$. In particular, the $\G^\infty$-regular embedded distributions
are precisely the smooth functions.

The distributional wave front set can be described (in fact, was
originally defined in \cite{Hoermander:71}) in terms of
characteristic sets of pseudodifferential operators: If $m\in\R$
arbitrary then for any $v \in \D'$
\beq \label{classical_WF_char}
  \WF(v) = \bigcap_{\substack{A \in \Psi^m\\ Av \in \Cinf}} \Char(A).
\eeq
Here, $\Char(A)$ is the \emph{characteristic set} of $A = a(x,D)$ ($a
\in S^m$ the space of smooth symbols of order $m$ and type $(0,1)$), i.e., the complement of those points $(x_1,\xi_1)$, where
$a$ is (micro-)elliptic in the sense that an estimate
$$
  |a(x,\xi)| \geq C (1 + |\xi|)^m
$$
holds for $x$ near $x_1$ and $\xi$ with $|\xi| \geq 1$ in some conic neighborhood of
$\xi_1$. Equation (\ref{classical_WF_char}) is an efficient distribution
theoretic lever when applied to questions of noncharacteristic
regularity or propagation of singularities for solutions to partial
differential equations, or to directly deduce geometric invariance
properties of the wave front set.

\paragraph{Slow scale micro-ellipticity:}

A key property for generalized symbols is that of slow scale
growth, which has already been used in several contexts 
of regularity theory (cf. \cite{HO:03, HOP:05, GGO:03, GH:05}).
The basic measurement is done by estimating with elements of the
following set of \emph{strongly positive slow scale nets}:
\[
  \Pi_\ssc := \{ (\omega_\eps)_\eps\in\R^{(0,1]}\, \colon
  \exists c > 0\, \forall p \ge 0\, \exists c_p > 0\,  \forall \eps : 
  c \le \omega_\eps \text{ and } \omega_\eps^p \le  c_p\, \eps^{-1} \}.
\]                      

\begin{definition}\label{def_ssc_symbols} 
 Let $m$ be a real number.
 The set $\Sscu^m(\Om\times\R^n)$ 
of slow scale nets of symbols of order $m$ consists of all 
$(a_\eps)_\eps\in {\mathcal{S}}^m[\Om\times\R^n]$ such that:
\begin{multline*}
 \forall K \Subset \Omega\ \exists (\omega_\eps)_\eps \in \Pi_\ssc
 \forall \alpha, \beta \in \N^n\, \exists c>0\,  \forall \eps:\\
 \sup_{x\in K,\xi\in\R^n}\lara{\xi}^{-m+|\alpha|}|\partial^\alpha_\xi\partial^\beta_x a_\eps(x,\xi)|\le c\, \omega_\eps.
\end{multline*}
The subset $\Nu^m(\Om\times\R^n)$ of negligible nets of symbols of order
$m$ is defined by $(a_\eps)_\eps\in
{\mathcal{S}}^m[\Om\times\R^n]$ with the property
\begin{multline*}
 \forall K \Subset \Omega\, \forall \alpha, \beta \in \N^n\,
 \forall q\in\N, \exists c>0\,  \forall \eps:\\
 \sup_{x\in K,\xi\in\R^n}\lara{\xi}^{-m+|\alpha|}|\partial^\alpha_\xi\partial^\beta_x a_\eps(x,\xi)|  \le c\, \eps^q.
\end{multline*}
The elements of the factor space $\Syscu^m(\Om\times\R^n) 
:= \Sscu^m(\Om\times\R^n) / \Nu^m(\Om\times\R^n)$
are the classes of  \emph{slow scale generalized symbols of order $m$}.
\end{definition} 

Slow scale generalized symbols enable us to design a simple and 
useful notion of micro-ellipticity (cf.\ \cite{GH:05}).
\begin{definition}
\label{def_micro_ellipticity}
Let $a \in \Syscu^m(\Om\times\R^n)$ and $(x_0,\xi_0) \in \CO{\Om}$. We
say that $a$ is slow scale micro-elliptic at $(x_0,\xi_0)$ if it has a
representative $(a_\eps)_\eps$ satisfying the following: there is a
relatively compact open neighborhood $U$ of $x_0$, a conic
neighborhood $\Gamma$ of $\xi_0$, and $(r_\eps)_\eps , (s_\eps)_\eps$
in $\Pi_\ssc$ such that 
\beq \label{estimate_below} 
  | a_\eps(x,\xi)|  \ge \frac{1}{s_\eps} \lara{\xi}^m\qquad\qquad (x,\xi)\in
  U\times\Gamma,\, |\xi| \ge r_\eps,\, \eps \in (0,1].  
\eeq
We denote by $\Ellsc(a)$ the set of all $(x_0,\xi_0) \in \CO{\Om}$
where $a$ is slow scale micro-elliptic. The symbol $a$ is called 
\emph{slow scale elliptic} if there exists $(a_\eps)_\eps \in\ a$ 
such that \eqref{estimate_below} holds at all points in $\CO{\Om}$.
\end{definition}
In case of classical symbols the notion of slow scale
micro-ellipticity coincides with the classical one, which is defined 
equivalently as the set of noncharacteristic points. 
More general definitions of ellipticity have been investigated in 
\cite{Garetto:04, GGO:03, HO:03, HOP:05}. But the above slow
scale variant has a stability property, similar to classical ellipticity. 
In fact, due to the overall slow scale
conditions in Definition \ref{def_ssc_symbols}, any symbol which is
slow scale micro-elliptic at $(x_0,\xi_0)$ fulfills the stronger
hypoellipticity estimates \cite[Definition 6.1]{GGO:03}; furthermore,
\eqref{estimate_below} is stable under lower order (slow scale)
perturbations (cf.\ \cite[Proposition 1.3]{GH:05}), and the 
simple slow scale ellipticity
condition in Definition \ref{def_micro_ellipticity} already guarantees
the existence of a parametrix. 

\paragraph{Pseudodifferential characterization of the generalized wave
  front 
set:}

Let $\Oprop{m}(\Om)$ denote the set of all properly
supported operators $a(x,D)$ where $a$ belongs to
$\Syscu^m(\Om\times\R^n)$. 
\begin{theorem}[\cite{GH:05}]
\label{theorem_wave_front}
Let $m\in\R$ arbitrary. For every $u\in\G(\Om)$
$$
  \WFg(u) = \hskip-5pt
 \bigcap_{\substack{a(x,D)\in\,\Oprop{m}(\Om)\\[0.1cm] a(x,D)u\, \in\,
     \Ginf(\Om)}}\hskip-5pt \compl{\Ellsc(a)}
 = \bigcap_{\substack{A\in\Psi^m(\Om)\\[0.1cm] Au\, \in\,
     \Ginf(\Om)}}\hskip-5pt  \Char(A).
$$  
\end{theorem}

\paragraph{Noncharacteristic $\Ginf$-regularity:} A  
first application of Theorem \ref{theorem_wave_front} is
the following extension of classical noncharacteristic 
regularity for solutions of pseudodifferential 
equations (cf.\ \cite[Theorem 18.1.28]{Hoermander:V3}). 
\begin{theorem}[\cite{GH:05}] If $P = p(x,D)$ is a properly supported 
pseudodifferential operator with slow scale symbol and $u\in\G(\Om)$ 
then 
\beq\label{nonchar_theorem}
  \WFg(P u) \subseteq \WFg(u) \subseteq
    \WFg(P u) \cup \compl{\Ellsc(p)}.
\eeq
\end{theorem}

As can be seen from various examples in \cite{HO:03},
relation (\ref{nonchar_theorem}) does not hold in general for regular symbols
$p$ which satisfy estimate (\ref{estimate_below}). In this sense, the overall
slow scale properties of the symbol are crucial in the above statement and are
not just technical convenience. In fact, adapting the reasoning in
\cite[Example 4.6]{HO:03} to the symbol $p_\eps(x,\xi) = 1 + c_\eps x^2$,
$c_\eps \geq 0$, we obtain the following: $p_\eps(x,\xi) \geq 1$ whereas the
unique solution $u$ of $p u = 1$ is $\Ginf$ if and only if $(c_\eps)_\eps$ is a
slow scale net.

\subsection{Propagation of $\Ginf$-singularities} 

\paragraph{1. The case of smooth coefficients:}
In \cite{GH:05} the pseudodifferential approach to wave front sets provides
the key tool to extend a distributional result on propagation of
singularities.  We allow for  
Colombeau generalized functions as solutions and initial values in first-order
strictly hyperbolic partial differential equations with smooth coefficients.
Hyperbolicity will be assumed with respect to time direction,
corresponding to the variable $t$. 
\begin{theorem}[\cite{GH:05}] \label{smooth_prop_thm}
Let $P(t,x,D_x)$ be a first-order partial differential operator with 
real principal symbol $p_1$ and coefficients in $\Cinf(\R\times\Om)$, which
are constant outside some compact subset of $\Om$.  
Let $u\in\G(\Om\times\R)$ be the (unique) solution to the homogeneous 
Cauchy problem
 \begin{align}
 D_t u + P(t,x,D_x) u &= 0 
                \label{first_order_equ}\\
 u(.,0) &= g \in\G(\Om) \label{initial_cond}.
 \end{align}
\begin{enumerate}
\item If $\Phi_t$ denotes the Hamilton flow corresponding to $p_1(t,.,.)$ on
$T^*(\Om)$ then we have for all $t\in\R$  
 \beq\label{WFg_flow}
    \WFg(u(.,t)) = \Phi_t\big( \WFg(g) \big).
 \eeq
\item Denote by $\ga(x_0,\xi_0)$ the maximal bicharacteristic 
curve (of $Q := D_t + P$) passing through
$(x_0,0;\xi_0,-P_1(0,x_0,\xi_0)) \in\Char(Q) \subseteq
\CO{\Om\times\R}$. 
Then the generalized wave front of $u$ is given by
 \beq
    \WFg(u) = \bigcup_{(x_0,\xi_0)\in\WFg(g)} \ga(x_0,\xi_0).
 \eeq
\end{enumerate}
\end{theorem}

\paragraph{2. Colombeau-type coefficients and generalized (null)
  bicharacteristics:}\ \\ 
We study global problems with $\Om = \R^n$. 
Consider the first-order partial differential operator 
\[
   P(t,x;D_x) = \sum_{j=1}^n a_j(x,t) D_{x_j} + a_0(x,t),
\]
where $a_k \in \G(\R^{n+1})$ ($k=0,\ldots,n$) and $a :=
(a_1,\ldots,a_n) \in \G(\R^{n+1})^n$ is real, in the sense that
each component has real-valued representatives. In addition, we make
the following assumptions on the coefficients:
\begin{enumerate}
\item $a_0,\ldots, a_n$ are equal to real constants for large $|x|$, 
\item $\d_{x_j} a_k$ ($k = 1,\ldots,n$; $j=1,\ldots,n$) as well as
  $a_0$ are of $\log$-type (i.e., satisfy asymptotic $L^\infty$-estimates
  with bound $O(\log(1/\eps))$).
\end{enumerate}
Then we deduce from \cite{LO:91} that the corresponding Cauchy problem
\begin{align}
  Qu := D_t u + P u & = f \in \G(\R^{n+1}) \label{gen_PDE}\\
     u\mid_{t=0} &= g \in \G(\R^n) \label{gen_initial_cond}
\end{align}
has a unique solution $u\in\G(\R^{n+1})$. The, now,
generalized principal symbols $p_1(t,.,.) = a(.,t) \cdot \xi 
\in \G(T^*(\Om))$ ($\forall t\in\R$) and 
$q_1 = \tau + p_1 \in \G(T^*(\R^{n+1}))$ define
the generalized Hamilton vector fields
\[
    H_{p_1(t)} = (a(.,t);-d_x a(.,t)\cdot \xi) \text{ and } 
    H_{q_1} = (a,1;-d_x a \cdot \xi,-\d_t a \cdot \xi)
\]
on $T^*(\R^n)$ and $T^*(\R^{n+1})$ respectively. Let
$(x_0,t_0;\xi_0,0)\in T^*(\R^{n+1})$ satisfy
$q_1(x_0,0,\xi_0,\tau_0) = 0$ (in $\widetilde{\C}$) and consider the
bicharacteristic system of generalized ordinary differential equations. Note that, due to the
simple symbol structure with respect to $\tau$, the $t$-component 
may be used as parameter along the curves; thus, the equations read
\begin{align}
  \dot{x}(t)   &= a(x(t),t) & x(0) &= x_0 \label{x_equ}\\
  \dot{\xi}(t) &= -d_x a(x(t),t) \cdot \xi(t) & \xi(0) 
                                                 &= \xi_0 \label{xi_equ}\\
  \dot{\tau}(t) &= -\d_t a(x(t),t) \cdot \xi(t)  & \tau(0) 
                                                 &= \tau_0.\label{tau_equ}
\end{align}
Here, (\ref{x_equ}) is independent of the other equations and is uniquely
solvable for $x \in \G(\R)^n$ by \cite[Theorem 1.5.2]{GKOS:01}. 
Since $d_x
a$ is of $\log$-type and Equation (\ref{xi_equ}) is linear with
respect to $\xi$, the standard Gronwall estimates (in terms of vector
and matrix norms) yield unique existence of the generalized solution 
$\xi \in \G(\R)^n$. Finally, $\tau$ is obtained by simply integrating
(\ref{tau_equ}). We summarize this in the following statement.
\begin{proposition} \label{gen_bichar}
Under Assumptions (i) and (ii) above for the coefficients of the 
operator $Q = D_t + P(t,x,D_x)$ there is a unique generalized 
bicharacteristic curve $\beta\in\G(\R)^{2n + 2}$, $\beta(t) 
= (x(t),t;\xi(t),\tau(t))$,
passing through every point $(x_0,0;\xi_0,\tau_0) \in T^*(\R^{n+1})$ at
which the principal symbol $q_1$ vanishes.
\end{proposition} 

\begin{remark} (i) As a matter of fact, (\ref{x_equ}) establishes a generalized
(c-bounded) flow on $\R^n$ in the sense of \cite{KOSV:04}, whereas
(\ref{x_equ}-\ref{xi_equ}) defines a (not necessarily c-bounded) flow
on $T^*(\R^n)$ 
as in \cite[Equations (1.33-35)]{GKOS:01} (and without
satisfying the hypothesis of Theorem 1.5.12 there, since the 
$(x,\xi)$-gradient of the vector field $H_{p_1(t)}$ need
not be of $\log$-type).\\
(ii) Observe that, by construction, we have that
\beq\label{Hamilton_null}
   q_1(x(t),t,\xi(t),\tau(t)) = 0 \text{ in } \widetilde{\R} \qquad
   \forall t\in\R,
\eeq
i.e., the Hamilton function vanishes along the generalized null bicharacteristics.
\end{remark}

We may then ask, if the generalized bicharacteristic flow still
describes the propagation of microlocal singularities in a similar way
as Theorem \ref{smooth_prop_thm} illustrates this for the case of smooth
coefficients. In other words: How does the generalized wave front set, 
$\WFg(u)$, of the solution $u$ to (\ref{gen_PDE}-\ref{gen_initial_cond}) 
relate to the generalized bicharacteristic curves (or the flow)
obtained in Proposition \ref{gen_bichar}? 

On a conceptual level, one immediately realizes that a direct set
theoretic relation is not meaningful since $\WFg(u) \subseteq T^*(\R^{n+1})
\isom \R^{2n+2}$, whereas the bicharacteristic curves give generalized
point values $\beta(t) \in \widetilde{\R}^{2n+2}$ at fixed $t\in\R$.
Here, we will only touch upon this new research issue by discussing a simple
example in space dimension $n = 1$.

\paragraph{The Hurd-Sattinger Example revisited:} Let
$\rho\in\Cinfc(\R)$ be symmetric, nonnegative, with $\supp(\rho)
\subseteq [-1,1]$ and $\int \rho =
1$. We use the (slow) scale $\ga_\eps := \log(1/\eps)$, put
$\rho^\eps(x) := \ga_\eps \rho(\ga_\eps x)$, and define ($H$ denoting the
Heaviside function)
\[
   \Theta := \cl{(H(-.) * \rho^\eps)_{\eps\in (0,1]}} \approx H(-x) 
\]
considered as element in $\G(\R^2)$ (i.e., independent of the
$t$-variable). Note that $\Theta$ as well as $\d_x \Theta = \Theta'$ are of
$\log$-type and that $\Theta$ equals $1$ when $x < -1$ and equals $0$ when
$x > 1$. 

Let $u\in\G(\R^2)$ be the solution to the Cauchy problem
\begin{align}
  Qu := D_t u + \Theta \d_x u + \Theta' u & = 0 \label{HS_PDE}\\
     u\mid_{t=0} &= {d_{-s_0}}  \label{HS_initial_cond}
\end{align}
where $s_0 > 0$ and $d_{-s_0} = \cl{(\rho^\eps(.+s_0))_\eps} \approx
\delta_{-s_0}$ (the Dirac measure at $-s_0$). The bicharacteristic
differential equations read
\begin{align}
  \dot{x}(t)   &= \Theta(x(t)) & x(0) &= x_0 \label{HS_x_equ}\\
  \dot{\xi}(t) &= -\Theta'(x(t)) \xi(t) & \xi(0) 
                                                 &= \xi_0 \label{HS_xi_equ}\\
  \dot{\tau}(t) &= 0  \quad \Rightarrow \quad \tau(t) = \tau_0 \quad \forall t.
\end{align} 
Integrating (\ref{HS_xi_equ}) we get
 \beq\label{xi_expr}
      \xi_\eps(t) = \xi_0 \exp\big(- \!\!
      \int\limits_0^t \!\Theta'_\eps(x_\eps(r))\, dr\big)
      =   \xi_0 \exp\big( \int\limits_0^t \!\rho_\eps(-x_\eps(r))\, dr\big) .
 \eeq
Note that, by construction, we have the following characteristic
relation along the bicharacteristics: for all $t$ and $q$
 \beq\label{char_rel}
    (q_1(x(t),t,\xi(t),\tau(t)))_\eps = \tau_\eps(t) +
    \Theta_\eps(x_\eps(t))\, \xi_\eps(t) = O(\eps^q) \quad (\eps \to 0).
 \eeq
A detailed description of the solution to (\ref{HS_x_equ}), the
characteristic flow, and an analysis of the generalized wave front set
of $u$ can be found in \cite{HdH:01}. To be
precise, the microlocal regularity properties are studied there on the
scale $\ga_\eps$, indicated by the notation $\WFg^\ga(u)$, in order to
correctly trace back the influence of the ``original distributional'' 
coefficient singularity at $x = 0$. Accordingly, the initial value
singularity at $x_0 = - s_0 < 0$ propagates towards the axis $x = 0$ with
speed $1$ as long as $x < 0$ and ``gets stuck'' at $x = 0$ from time $t =
s_0$ on. The wave front set turns out to be
\begin{multline}\label{WF_sol}
  \WFg^\ga(u) = \big(\{ t = x + s_0, x < 0\} \times \{ \xi + \tau = 0,
  (\xi,\tau) \not= (0,0)\}\big) \\ 
  \cup \big(\{ (0,s_0) \} \times (\R^2\setminus\{(0,0)\})\big) \\ 
  \cup \big(\{ x=0, t > s_0 \} \times \{ \tau = 0,  \xi \not= 0\}\big).
\end{multline}
We obtain the following picture for the representatives and the limit 
(as $\eps \to 0$) of the (bi)characteristic curve passing through
$(-s_0,0)$ (resp. $(-s_0,0;\xi_0,\tau_0)$):\\ 
\begin{minipage}{7cm}{\epsfig{file=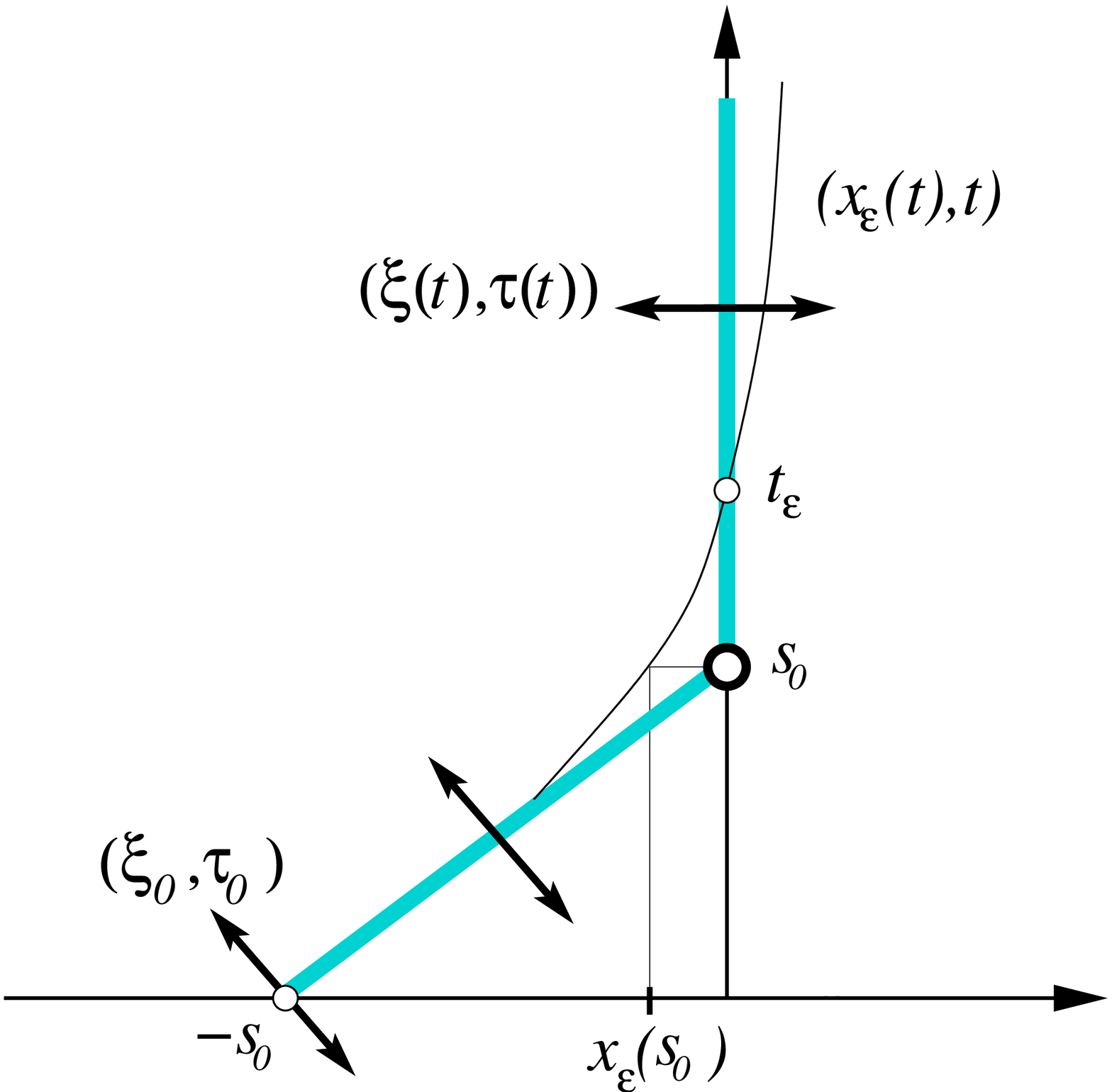, scale=0.3}}\end{minipage}
\begin{minipage}{5cm} Here, $t_\eps$ is unique with the
property $x_\eps(t_\eps) = 0$ and $t_\eps \to s_0$ as $\eps \to
0$, where the characteristic limit curve (thick gray line) is
non-differentiable.  

The admissible directions for $(\xi_0,\tau_0)$ are simply described by
$\tau_0 / \xi_0 = \pm 1$. Since  $\tau(t) = \tau_0$ for all $t$, it
remains to investigate $\lim_{\eps \to 0} \xi_\eps(t)$ for the cases
$t < s_0$, $t = s_0$, and $t > s_0$.
\end{minipage} 

As long as $t < s_0$ we directly obtain from 
equation (\ref{xi_expr}), the fact that $x_\eps\to t - s_0 < 0$ as $\eps\to 0$ (cf.\ \cite[Subsection 4.2]{HdH:01}), and
the support properties 
of $\rho$ that eventually (as soon as $\eps$ is sufficiently small)
$\xi_\eps(t) = \xi_0$, which is in accordance with the wave front set
in (\ref{WF_sol}). 

If $t > s_0$ then $\lim_{\eps \to 0} \dot{x}_\eps(t) = 0$ since the
characteristic limit curve is a vertical line in that region. Thus,
noting that $\dot{x}_\eps(t) \geq 0$, upon inserting (\ref{HS_x_equ})
into (\ref{char_rel}) we conclude that $\lim_{\eps \to 0} \xi_\eps(t)
= \pm \infty$, the sign corresponding to that of $-\tau_0$. This
yields the horizontal directions of the wave front set  (\ref{WF_sol})
as limit directions of $(\xi_\eps(t),\tau_\eps(t))$.

Finally, consider $t = s_0$, so we are right at the kink in the limit
of the characteristic curves. We note that
$\dot{x}_\eps(t_\eps) = \Theta_\eps(0) = \int_{-\infty}^0 \rho
= 1/2$ and, by monotonicity of $\Theta_\eps$ and $x_\eps$, the
inequality $\dot{x}_\eps(s_0) \geq  \dot{x}_\eps(t_\eps) \geq 0$ is
valid. Hence (\ref{char_rel}) and (\ref{xi_expr}) imply $|\tau_0| =
|\xi_\eps(s_0)| \dot{x}_\eps(s_0) \geq |\xi_\eps(s_0)| \dot{x}_\eps(t_\eps) =
|\xi_\eps(s_0)| / 2 \geq |\xi_0| / 2 > 0$, so that 
$$
    \tau_\eps(s_0) = \tau_0, \quad 
    |\xi_0| \leq |\xi_\eps(s_0)| \leq 2 |\tau_0|   
    \qquad \forall \eps \in   (0,1], 
$$
which shows that at the base point $(0,s_0)$ the possible
accumulation points of the bicharacteristic flow can cover only a true
subcone of all the directions actually occurring in the wave front set
(\ref{WF_sol}) at this singularity. In this sense, \emph{the generalized
bicharacteristic flow-out fails to generate the full wave front set}.



\bibliographystyle{abbrv}
\newcommand{\SortNoop}[1]{}

\end{document}